\documentclass{amsproc}
\usepackage{color}
\newtheorem{Theorem}{Theorem}
\newtheorem{theorem}{Theorem}[section]
\newtheorem{lemma}[theorem]{Lemma}

\newtheorem{corollary}[theorem]{Corollary}

\theoremstyle{definition}
\newtheorem{definition}[theorem]{Definition}

\newtheorem{notation}[theorem]{Notation}

\theoremstyle{remark}

\numberwithin{equation}{section}

\begin{document}

\title[Invariable generation]{Invariable generation of certain groups of piecewise linear
homeomorphisms of the interval}


\author{Yoshifumi Matsuda}
\address{Department of Physics and Mathematics, College of Science and
Engineering, Aoyama Gakuin University, 5-10-1 Fuchinobe, Chuo-ku,
Sagamihara-shi, Kanagawa, 252-5258}
\email{ymatsuda@gem.aoyama.ac.jp}
\author{Shigenori Matsumoto}
\address{Department of Mathematics, College of
Science and Technology, Nihon University, 1-8-14 Kanda-Surugadai,
Chiyoda-ku, Tokyo, 101-8308 Japan}
\email{matsumo@math.cst.nihon-u.ac.jp}

\thanks{The first author is partially supported by Grant-in-Aid for
{{Young Scientists (B)}}  No.\ 25800036, and the second by {{Grant-in-Aid for
Scientific Research}} (C) No.\ 25400096.}

\date{\today}

\newcommand{\AAA}{{\mathcal A}}
\newcommand{\BBB}{{\mathbb B}}
\newcommand{\LL}{{\mathcal L}}
\newcommand{\MCG}{{\rm MCG}}
\newcommand{\PSL}{{\rm PSL}}
\newcommand{\R}{{\mathbb R}}
\newcommand{\Z}{{\mathbb Z}}
\newcommand{\XX}{{\mathcal X}}
\newcommand{\per}{{\rm per}}
\newcommand{\N}{{\mathbb N}}

\newcommand{\PP}{{\mathcal P}}
\newcommand{\GG}{{\mathbb G}}
\newcommand{\FF}{{\mathcal F}}
\newcommand{\EE}{{\mathbb E}}
\newcommand{\BB}{{\mathbb B}}
\newcommand{\CC}{{\mathcal C}}
\newcommand{\HH}{{\mathcal H}}
\newcommand{\UU}{{\mathcal U}}
\newcommand{\oboundary}{{\mathbb S}^1_\infty}
\newcommand{\Q}{{\mathbb Q}}
\newcommand{\DD}{{\mathcal D}}
\newcommand{\rot}{{\rm rot}}
\newcommand{\Cl}{{\rm Cl}}
\newcommand{\Index}{{\rm Index}}
\newcommand{\Int}{{\rm Int}}
\newcommand{\Fr}{{\rm Fr}}
\newcommand{\ZZ}{\Z[2^{-1}]}
\newcommand{\II}{{\mathcal I}}
\newcommand{\JJ}{{\mathcal J}}
\newcommand{\g}{{\rm gen}}
\newcommand{\KK}{{\mathcal K}}
\newcommand{\OO}{{\mathcal O}}

\date{\today }

\maketitle

\begin{abstract}
Let $P$ be the group of all the orientation preserving piecewise linear
 homeomorphisms of the interval $[0,1]$.
Given any $a>1$, let $P^a$ be the subgroup of $P$
consisting of all the elements with slopes in $a^\Z$,
and let  $P^\Q$ be the subgroup of $P$ consisting of all the elements
 with slopes and breaks in $\Q$.
We show that the groups $P$, $P^a$, $P^\Q$, as well as Thompson group $F$,
are invariably generated.
\end{abstract}

\section{Introduction}  
The concept of invariable generation for a group $G$ was introduced by
J. Wiegold in \cite{Wiegold1}.

\begin{notation}\label{n1}\footnote{The convention $g^h$ is not the same as the customary one.}
For elements $g$ and $h$ of a group $G$, and a subgroup $H$ of $G$, we
 denote:
$$
h^g=ghg^{-1}, \ \ h^G=\{h^g\mid g\in G\}, \ \ H^g=\{h^g\mid h\in H\}.
$$
\end{notation}

\begin{definition}
 (1) A subgroup $H$ of $G$ is called {\em classful} if 
$H\cap  g^G\neq\emptyset$ for any $g\in G$, or equivalently, 
\begin{equation}\label{e1}
 \bigcup_{g\in
 G}H^g=G.
\end{equation}

(2) A group $G$ is said to be {\em invariably generated} if there {{are}} no
classful subgroups other than $G$ itself.
\end{definition}

Any finite group is invariably generated, as is shown by a counting
argument on (\ref{e1}). Much easier is the fact that
 any abelian group is invariably
generated. In \cite{Wiegold1}, it is shown that the invariable
generation is extension closed.
Therefore any virtually solvable group is invariably generated.  It is
also projection closed.
Given a prime number $p>10^{75}$, an infinite group whose arbitrary
proper nontrivial subgroup is of order $p$ is constructed in \cite{O}.
Such groups are necessarily generated by arbitrary two elements
not from the same proper subgroup, and is invariably
generated, provided there are more than one nontrivial conjugacy
classes. The Grigorchuk group \cite{G} is also invariably generated
\cite{KLS}.

However the invariable generation is not subgroup closed: an example is given in
\cite{Wiegold2}.
It is also not direct union closed: the group of the permutations
of $\N$ with finite support is not invariably generated, 
since the stabilizer of $1\in\N$ is classful.
Infinite groups with one nontrivial conjugacy class, constructed in
\cite{HNN}, {{are}} not invariably generated.
Free groups of generators $\geq2$ are not invariably generated
\cite{Wiegold1}. More generally, {{nonelementary convergence groups}} 
are not invariably generated \cite{Gelander}.
{{Acylindrically hyperbolic groups are not invariably generated \cite{BF}.}} 
Invariable generation of linear groups are discussed in \cite{KLS}. 

The current paper is concerned with groups of
piecewise linear (PL) homeomorphisms of the interval. 

\begin{Theorem}\label{t3}
 Thompson group $F$ is invariably generated.
\end{Theorem}

Our method cannot give the finite
invariance generation of $F$ obtained in \cite{GGJ}.
Denote by $P$ the group formed by all the orientation preserving PL homeomorphisms of the unit interval $[0,1]$, and by $P^\Q$ the
subgroup of $P$ formed by elements with slopes and breaks in $\Q$.
Fix an arbitrary
real number $a>1$. Let
$P^a$ be the subgroup of $P$ consisting of all the elements with slopes
in $a^\Z$.

\begin{Theorem}\label{t1} 
 The group $P^a$ is invariably generated.
\end{Theorem}

\begin{Theorem}\label{t2}
 The group $P$ is invariably generated.
\end{Theorem}

\begin{Theorem}\label{t4}
 The group $P^\Q$ is invariably generated.
\end{Theorem}

The proofs of the above theorems are quite similar. In Section 2, we
summarize conditions for a subgroup $G$ of $P$ to be invariably generated.
In later sections we show that $F$, $P^a$, $P$ and $P^\Q$ satisfy
these conditions independently.

\section{Conditions for invariable generation}

Let $G$ be any subgroup of the group $P$ of all the orientation preserving
PL homeomorphisms of the interval $[0,1]$. We shall raise three conditions for $G$
to be invariably generated. Let $X$ be a dense subset of $(0,1)$ which
is left invariant by $G$, and let $X^*=X\cup\{0,1\}$.
A closed interval $I\subset[0,1]$ is called an {\em $X$-interval} (resp.\ {\em
$X^*$-interval}) if the endpoints of $I$ are contained in $X$ (resp.\ $X^*$).

\begin{definition}\label{d1}
 For an $X^*$-interval $I$, let us denote
$$
G(I)=\{g\vert_I\mid g\in G,\ {\rm Supp}(g)\subset I\}\ \mbox{ and }\
G\vert_I=\{g\vert_I\mid g\in G,\ g(I)=I\}.
$$
\end{definition}

The first condition is to fix the relation between $G$ and $X$.

\bigskip
\bf Condition A: \rm
(1) The breaks of any $g\in G$ are contained in $X$.

(2) The group $G$ acts on $X$ transitively.

(3) For any $X$-interval $I$, $G\vert_I=G(I)$.

(4) For any $X^*$-interval $I$, there is a PL homeomorphism
$\psi_I:[0,1]\to I$ such that $\psi_I(X^*)=X^*\cap I$ and $G^{\psi_I}=G(I)$.

\bigskip
The other two conditions are concerned with an arbitrary classful
subgroup $H$ of $G$.

\bigskip
\bf Condition B: \rm Any classful subgroup $H$ acts on $X$ transitively.

\begin{definition}\label{d2}
 For an $X^*$-interval $I$, let us denote
$$
H(I)=\{h\vert_I\mid h\in H,\ {\rm Supp}(h)\subset I\}\ \mbox{ and }\
H\vert_I=\{h\vert_I\mid h\in H,\ h(I)=I\}.
$$
\end{definition}

\bf Condition C: \rm For any classful subgroup $H$, there is an $X$-interval
$I_0$ such that $H\vert_{I_0}=G(I_0)$.

\bigskip
In this section, we show that if a subgroup $G$ of $P$ satisfies
conditions A, B and C, then $G$ is invariably generated. 
\em Henceforth in this section, we assume $G$ satisfies conditions A, B
and C.
\rm  For $f\in G$, define
$s(f)\in [0,1]$ by $$s(f)=\sup\{s\mid f\vert_{[0,s]}={\rm id}\}.$$
By condition A(1), $s(f)$ is contained in $X^*$.
Let $H$ be an arbitrary classful subgroup of $G$

\begin{lemma}\label{l8}
For any $X^*$-interval $I=[t,1]\subset[0,1]$,
 $H(I)$
is a classful subgroup of $G(I)$.
\end{lemma}

\sc Proof: \rm Given any $f\vert_{I}\in G(I)$ where $f\in G$ with ${\rm Supp}(f)\subset I$, let us show that there is $g\in G$ such that ${\rm Supp}(g)\subset I$ and $f^{g}\in H$.
Since $H$ is 
classful, there is $g_1\in G$ such that $f^{g_1}\in H$. Now 
$s(f^{g_1})=g_1(s(f))\in X$. Notice that $s(f)\geq t$ since 
${\rm Supp}(f)\subset I$.
 By condition B, there is $g_2\in H$ such that
$g_2(g_1(s(f)))=s(f)$. Then $f^{g_2g_1}=(f^{g_1})^{g_2}$ is an
element in $H$ since $f^{g_1}\in H$ and $g_2\in H$.
Moreover $f^{g_2g_1}\vert_{I}$ belongs to
 $H(I)$ since it satisfies $s(f^{g_2g_1})=s(f)\geq t$. 
Notice that  $s(f)$ is a fixed point of $g_2g_1$.
Now by condition A(3), there is an element $g\in G$ which is
 the identity on $[0,s(f)]$ and is equal to
$g_2g_1$ on $[s(f),1]$. 
 Then we have $f^{g}=f^{g_2g_1}\in H$,
as is required. \qed

\bigskip
For any $X^*$-interval $I=[t,1]$, choose a
PL homeomorphism $\psi_I:[0,1]\to I$
 such that $\psi_I(X)=X\cap (t,1)$ and $G^{\psi_I}=G(I)$ (condition A(4)).
 Since $H(I)$ is classful in $G(I)$,
 $H(I)^{\psi^{-1}}$ is 
 classful in $G(I)^{\psi^{-1}}=G$. By condition B, $H(I)^{\psi^{-1}}$ acts 
transitively on $X$. Therefore $H(I)$ acts transitively on $X\cap(t,1)$.
This way we get the following lemma.

\begin{lemma}\label{l9}
 The classful subgroup $H$ acts doubly transitively on $X$.
\qed
\end{lemma}

By the same argument as Lemma \ref{l8},
 applied to the inclusion of an
$X$-interval $I=[t,t']$
into an $X^*$-interval $=[t,1]$, we get the following.

\begin{lemma}\label{l10}
 For any $X$-interval $I$ of $(0,1)$, the group $H(I)$ is
a classful subgroup of $G(I)$. \qed
\end{lemma}

We shall discuss consequences of condition C. Let $I_0$ be an
 $X$-interval 
 such that $H\vert_{I_0}=G(I_0)$.

\begin{lemma}\label{l12}
 We have
$H(I_0)=G(I_0)$.
\end{lemma}

\sc Proof: \rm Choose an arbitrary $f\vert_{I_0}\in G(I_0)$, where $f\in G$ with ${\rm Supp}(f)\subset I_0$. 
Then since $H(I_0)$ is classful in $G(I_0)$ (Lemma \ref{l10}),
there is $g\in G$ such that ${\rm Supp}(g)\subset I_0$ and $f^g\in H$. Since 
$g\vert_{I_0}\in G(I_0)=H\vert_{I_0}$,
there is $h\in H$ such that 
$h(I_0)= I_0$
and $h\vert_{I_0}=g\vert_{I_0}$. Then since ${\rm Supp}(f)\subset I_0$,
we have $f^g=f^h$, and hence 
$f=(f^g)^{h^{-1}}\in H$. This, together with the assumption
 ${\rm Supp}(f)\subset I_0$, implies that $f\vert_{I_0}\in H(I_0)$. \qed

\begin{corollary}\label{c3}
 For any $X$-interval $I$ in $(0,1)$, we have $H(I)=G(I)$.
\end{corollary}

\sc Proof. \rm By double transitivity of the action of $H$ on $X$
(Lemma \ref{l9}), there is $h\in H$ such that $h(I_0)=I$.
Now 
$$
H(I)=H(I_0)^h=G(I_0)^h=G(I),$$
as is required. \qed

\bigskip
Finally we shall prove that $H=G$.
Let 
$$G_0=\{g\in G\mid g'(0)=g'(1)=1\},\ \ 
H_0=\{h\in H\mid h'(0)=h'(1)=1\}.
$$

Let $\{J_n\}_{n\in\N}$ be an increasing sequence of $X$-intervals such that 
$\cup_nJ_n=(0,1)$. We have
$$
H_0=\bigcup_{n\in\N}H(J_n)\ \mbox{ and }\
G_0=\bigcup_{n\in\N}G(J_n).$$
Since by the previous lemma, $H(J_n)=G(J_n)$ for any $n\in\N$, 
we get $G_0=H_0$.

Now for any $f\in G$, there is $g\in G$ such that $f^g\in H$. But then
$f^gf^{-1}=[g,f]\in G_0\subset H$, and therefore $f\in H$.
This finishes the proof that if $G$ satisfies conditions A, B and C,
then $G$ is invariably generated.

\bigskip
In the rest of the paper, we use the following terminology.
\begin{definition}\label{elz}
 For $f\in P$, an interval $[0,\epsilon]$ or $[1-\epsilon,1]$ on which
$f$ is linear is called an {\em end linear zone} of $f$.
\end{definition}.

\section{The group $F$}

Let us denote by $\ZZ\subset\R$ the set of dyadic rationals.
Thompson group $F$ is the subgroup of $P$ consisting of all the
elements with slopes in $2^\Z$ and
breaks in $\ZZ$. For $F$, we define $X$ in the previous section
as $X=\ZZ\cap(0,1)$.
It is well known that $F$ satisfies condition A. For A(4), we can
take $\psi_I$ to be any PL homeomorphism from $[0,1]$ to $I$ with slopes
 in $2^\Z$ and breaks in $\ZZ$. See \cite{Belk} for the existence.

\begin{definition}
 \label{d2}
Define a homomorphisms $\alpha:F\to\Z^2$  by
$$\alpha(f)=(\log_2f'(0),\log_2f'(1)).$$
\end{definition}
Notice that $f\in{\rm Ker}(\alpha)$ if and only if
${\rm Supp}(f)\subset(0,1)$.
It is well known \cite{Belk} that ${\rm Ker}(\alpha)=[F,F]$.
Of course $\alpha$ is a class function: $\alpha(f)=\alpha(f^g)$.
Let 
$$F_{1,-1}=\{f\in F\mid f(x)>x,\forall x\in(0,1),\ \alpha(f)=(1,-1)\}.$$
Given $f\in F_{1,-1}$, 
 points $2^{-i}$ from an end linear zone of $f$ at $0$ are contained in a
 single orbit of the $\langle f\rangle$-action.
Their images  by high iterates of $f$ which lie in an end linear zone
at $1$ are of the form $1-k2^{-j}$ for some positive odd integer $k$.

\begin{definition}
 Define a map $\beta:F_{1,-1}\to 2\N-1$, 
by setting $\beta(f)$ to be the above odd integer $k$.
\end{definition}

\begin{lemma}\label{l3}
 The map $\beta$ is surjective. 
\end{lemma}

\sc Proof. \rm Let $k\in2\N-1$ be given.
Choose a large integer $j$, and define $g\in F$ by setting
$$ g(x)=2x \ \ \mbox{on} \ \ [0,2^{-j}], $$
$$ g(x)=2^{-1}(x-1)+1 \ \ \mbox{on} \ \ [1-k 2^{-j}.1], \ \ 
$$
and $g$ is a PL homeomorphism with slopes in $2^\Z$ and breaks in $\ZZ$
 from the interval
$[2^{-j},2^{-j+1}]$ to
$[2^{-j+1},1-k 2^{-j}]$.
Then we have $\beta(g)=k$. \qed

\begin{lemma}\label{l4}
 The map $\beta$ is class invariant, that is, $\beta(g^f)=\beta(g)$ for
any $g\in F_{1,-1}$ and any $f\in F$. 
\end{lemma}

\sc Proof. \rm Assume $\beta(g)=k\in2\N-1$ for $g\in F_{1,-1}$.
Then there is an orbit $\OO$ of $g$ which contains
$2^{-j}$ and $1-k2^{-j}$ for any large $j$. Choose an
arbitrary element $f\in F$ and assume that the slopes of $f$ are $2^{j_0}$ near $0$
and $2^{j_1}$ near $1$. Then $f$ maps $\OO$ to an orbit of $g^f$
which contains $2^{-j+j_0}$ and $1-k2^{-j+j_1}$ for any large
$j$, showing that $\beta(g^f)=k$. \qed

\bigskip
Let $H$ be an arbitrary classful subgroup of $F$.
\begin{corollary}\label{c1}
The map $\beta$ restricted to 
$H\cap  F_{1,-1}$ is surjective onto 
$2\N-1$. \qed
\end{corollary}

The next lemma shows that condition B of Section 2 is satisfied by $F$.
\begin{lemma}\label{l7}
 The classful subgroup $H$ acts transitively on $X=\ZZ\cap(0,1)$.
\end{lemma}

\sc Proof. \rm By Corollary \ref{c1}, there is
 an element $h_0\in H$ such that
$\beta(h_0)=1$. Thus for any large $j$, the points $2^{-j}$, as well as 
$1-2^{-j}$, are
on one orbit of $h_0$.
Again  by Corollary \ref{c1}, the $H$ orbit of these points
contains $1-k a^{-j}$ for any $k\in2\N-1$ and any large $j$. 
Applying negative iterates of $h_0$,
we get that the $H$ orbit contains all the points in $X$.
\qed

\bigskip
We need more in order to establish condition C for $F$.
For $n$ large, let $I_n=[2^{-n-1},2^{-n}]$, $J_n=[1-2^{-n},1-2^{-n-1}]$ and 
let $\phi_n:[0,1]\to I_n$, $\psi_n:[0,1]\to J_n$ be the orientation
preserving surjective linear map  of slope $2^{-n-1}$. 
Let
$$F_{1,-1,1}=\{g\in F_{1,-1}\mid\beta(g)=1\}.$$
Given any $g\in F_{1,-1,1}$, if we choose $n$ large
enough, some iterate
$g^N$ maps $I_n$ onto $J_n$. The map 
$\psi_n^{-1}\circ g^N\circ\phi_n$ is independent of the choice of $n$.
In fact, if $k>0$, $g^k\phi_{n+k}=\phi_n$ and $g^k\psi_n=\psi_{n+k}$.
Therefore we have
$$
\psi_{n+k}^{-1}g^{N+2k}\phi_{n+k}=(\psi_{n+k}^{-1}g^k)g^N(g^k\phi_{n+k})
=\psi_n^{-1}g^N\phi_n.$$
Notice also that $\psi_n^{-1}g^N\phi_n$ is an element of $F$.

\begin{definition}
Define a map
$\gamma:F_{1,-1,1}\to F$ by $\gamma(f)=\psi_n^{-1}\circ
g^N\circ\phi_n$.
\end{definition}

\begin{lemma}\label{l5}
The map $\gamma$ is surjective. 
\end{lemma}

We shall adopt a bit longer proof, which is applicable also to the group
$P^a$ in the next section. 

\bigskip
\sc Proof. \rm Choose an arbitrary element $g\in F_{1,-1,1}$ which 
is linear on $[0,2^{-n}]$ and $[1-2^{-n},1]$. There is $N>0$
such that $g^N$ maps $I_n$ onto $J_n$. Let $f_0=\psi_n^{-1}g^N\phi_n\in F$.
Any element of $F$ can be written as $ff_0$ for some $f\in F$.
The map $\hat f=\psi_n f\psi_n^{-1}$ is a PL homeomorphism of the
interval $J_n$ with slopes in $2^\Z$ and breaks in $\ZZ$. 
Define an element $g_1\in F_{1,-1,1}$ to be equal to 
$\hat fg$ on $J_{n-1}$ and equal to $g$ elsewhere. Notice that $\hat fg$ is
still linear on $[0,2^{-n}]$ and $[1-2^{-n},1]$. We also have 
$$
\psi_n^{-1}g_1^N\phi_n=\psi_n^{-1}\hat fg^N\phi_n
=(\psi_n^{-1}\hat f\psi_n)(\psi_n^{-1}g^N\phi_n)=ff_0.$$
Since $ff_0$ is an arbitrary element of $F$, we are done. \qed

\begin{lemma}\label{l6}
 The map $\gamma$ is class invariant. Precisely, if $g\in F_{1,-1,1}$
and $f\in F$, then $\gamma(g^{f})=\gamma(g)$.
\end{lemma}

\sc Proof. \rm Choose $n$ large enough so that
 $g$ and $f$ are linear on the intervals $[0,2^{-n}]$ and $[1-2^{-n},1]$.
Since $g\in F_{1,-1,1}$, some iterate $g^N$ of $g$ maps 
$I_n$ to $J_n$.
Put $k=\gamma(g)=\psi_n^{-1}g^N\phi_n$, and let us show that 
$\gamma(g^{f})=k$. We assume $f$ is of slope $2^{j_0}$ on $[0,2^{-n}]$
and of slope $2^{j_1}$ on $[1-2^{-n},1]$.
Then $g^{f}$ is linear (of slope $2$) on $[0,2^{-n+j_0}]$, and is
linear (of slope $2^{-1}$) on $[1-2^{-n+j_1},1]$. 
The map $g^{f}$ maps $I_{n-j_0}$ onto $J_{n-j_1}$.
Since $\phi_{n-j_0}=f\phi_n$ and $\psi_{n-j_1}=f\psi_n$, we have
$$
\psi_{n-j_1}^{-1}(g^{f})^{N}\phi_{n-j_0}=
\psi_n^{-1}f^{-1}(fgf^{-1})^Nf\phi_n=\psi_n^{-1}g^N\phi_n=k.$$
If $n$ is big enough compared with $j_0$ and $j_1$, we have
$\phi_{n-j_0}=(g^{f})^{j_0-j_1}\phi_{n-j_1}$. Therefore
$$k=\psi_{n-j_1}^{-1}(g^{f})^{N}\phi_{n-j_0}
=\psi_{n-j_1}^{-1}(g^{f})^{N+j_0-j_1}\phi_{n-j_1}.
$$
This shows $\gamma(g^{f})=k$, as is required. \qed

\bigskip

\begin{corollary}\label{c2}
 The map $\gamma$ restricted to $H_{1,-1,1}=H\cap F_{1.-1,1}$ is surjective onto
 $F$.
\qed
\end{corollary}

Fix once and for all an element $h_0\in H_{1,-1,1}$ such that $\gamma(h_0)={\rm id}$.
Thus there is $n>0$ such that $h_0$ is linear on $[0,2^{-n}]$ and 
$[1-2^{-n},1]$, that some iterate $h_0^{N}$ maps $I_n$ onto $J_n$ and that
$\psi_n^{-1}h_0^{N}\phi_n={\rm id}$.
The next lemma shows that the group $F$ satisfies condition C of
Section 2.

\begin{lemma}\label{l13} 
We have $H\vert_{I_n}=F(I_n)$.
\end{lemma}

\sc Proof: \rm Choose an arbitrary element $\hat f\in F(I_n)$ and let
$f=\phi_n^{-1}\hat f\phi_n\in F$.
By Corollary \ref{c2}, there is $h_1\in H_{1,-1,1}$ such that
$\gamma(h_1)= f$. More precisely, for some big $m>0$, there is $N>0$
such that $h_1^N(I_m)=J_m$ and that $\psi_m^{-1}h_1^N\phi_m= f$. One
can choose $m$ to be greater than $n$ in the lemma.
Then some iterate $h_0^{N_0}$ of $h_0$ maps $I_m$ onto
$J_m$ and $\psi_m^{-1}h_0^{N_0}\phi_m$ is still the identity.
Thus 
$$\phi_m^{-1}h_0^{-N_0}h_1^N\phi_m=
(\psi_m ^{-1}h_0^{N_0}\phi_m)^{-1}(\psi_m^{-1}h_1^N\phi_m)={\rm id}^{-1}f=f,$$
 and since
$h_0^{n-m}\phi_n=\phi_m$,
$$\phi_n^{-1}h_0^{m-n}(h_0^{-N_0}h_1^N)h_0^{n-m}\phi_n=f.$$
But this means 
$$h_0^{m-n}(h_0^{-N_0}h_1^N)h_0^{n-m}\vert_{I_n}=\hat f.$$
Since $\hat f\in F(I_n)$ is arbitrary and the LHS is in $H\vert_{I_n}$,
we are done.
\qed

\section{The group $P^a$}

Let $a>1$ be an arbitrary real number. 

\begin{definition} Given two compact intervals
$I$ and $J$, we denote by $PL^a(I,J)$ the space
of the PL homeomorphisms from $I$ to $J$ with slopes 
in $a^\Z$. Such a map is called a $PL^a$ {\em homeomorphism}.
\end{definition}

\begin{lemma}\label{l1}
For any compact interval $I$ and $J$, the space $PL^a(I,J)$ is nonempty.
\end{lemma}

\sc Proof. \rm Let $I=[p,q]$ and $J=[r,s]$. Consider
a line $L\subset\R^2$ of slope $a^n$, $n>1$, passing through the point $(p,r)$, and
another line $L'$ of slope $a^{-m}$, $m>1$, passing through
$(q,s)$. If $n$ and $m$ are sufficiently large, $L$ and $L'$ intersect at a point
in the open rectangle $(p,q)\times(r,s)$, yielding the graph of a desired map in
$PL^a(I,J)$. \qed

\begin{definition}
Define a group $P^a$ by  $P^a=PL^a([0,1],[0,1])$.
\end{definition}

We choose $X=(0,1)$ in condition A. Then the group $P^a$ satisfies 
A(1), A(2) and A(3) by virtue of Lemma \ref{l1}. For A(4), we just take
 $\psi_I:[0,1]\to I$ to be the orientation preserving linear homeomorphism.
Therefore in this section, $X^*$-intervals are 
just closed intervals. 
In the rest we shall establish conditions B and C for $P^a$ by almost
the same method as in Section 3.

Define a homomorphism $\alpha:P^a\to\Z^2$  by
$$\alpha(f)=(\log_af'(0),\log_af'(1)).$$
Clearly  $\alpha$ is a surjective class function.
Let 
$$P^a_{1,-1}=\{g\in P^a \mid g(x)>x,\forall x\in(0,1),\alpha(g)=(1,-1)\}.$$
Given $g\in P^a_{1,-1}$,  points $a^{-i}$ for $i$ large
are contained in a single orbit of $g$.
Consider their images by high iterates of $g$ which are near 1.
They are of the form $1-\xi a^{-j}$ for some number
$\xi\in(a^{-1},1]$.

\begin{definition}
  Define a map $\beta:P^a_{1,-1}\to (a^{-1},1]$, 
by setting $\beta(g)$ to be the above number $\xi$.
\end{definition}

Then one can show that 
 the map $\beta$ is a surjective class function just as Lemmas \ref{l3}
 and \ref{l4} in Section 3.
In particular, the map $\beta$ restricted to 
$H\cap  P^a_{1,-1}$ is surjective onto 
$(a^{-1},1]$, where $H$ is an arbitrary classful subgroup of $P^a$.
Then by the same method as Lemma \ref{l7}, we get the following lemma,
which establishes condition B.

\begin{lemma}
 Any classful subgroup $H$ acts transitively on $(0,1)$.
\end{lemma}

For a positive integer
$n$, let $I_n=[a^{-n-1},a^{-n}]$, $J_n=[1-a^{-n},1-a^{-n-1}]$ and 
let $\phi_n:[0,1]\to I_n$, $\psi_n:[0,1]\to J_n$ be the orientation
preserving linear homeomorphism {\em of the same slope $a^{-n}(1-a^{-1})$}. 
Let
$$P^a_{1,-1,1}=\{g\in P^a_{1,-1}\mid\beta(g)=1\}.$$
Given any $g\in P^a_{1,-1,1}$, if we choose $n$ large
enough, then $g$ is linear on the intervals $[0,a^{-n}]$ and
$[1-a^{-n},1]$. By the definition of $P^a_{1,-1,1}$, some iterate
$g^N$ sends $I_n$ to $J_n$, and the map 
$\psi_n^{-1}\circ g^N\circ\phi_n$ is independent of the choice of $n$.
Notice also that $\psi_n^{-1}g^N\phi_n$ is an element of $P^a$,
since $\phi_n$ and $\psi_n$ are linear
homeomorphisms of the same slope.

\begin{definition}
Define a map
$\gamma:P^a_{1,-1,1}\to P^a$ by $\gamma(f)=\psi_n^{-1}\circ
g^N\circ\phi_n$.
\end{definition}

One can show that the map $\gamma$ is a surjective class function
just as in Lemmas \ref{l5} and \ref{l6}.
Fix once and for all an element $h_0\in H\cap P^a_{1,-1,1}$ such that 
$\gamma(h_0)={\rm id}$.
Thus there is $n>0$ such that $h_0$ is linear on $[0,a^{-n}]$ and 
$[1-a^{-n},1]$, that some iterate $h_0^{N}$ maps $I_n$ onto $J_n$ and that
$\psi_n^{-1}h_0^{N}\phi_n={\rm id}$.
Just as in  Lemma \ref{l13}, we get the following lemma which establishes
condition C.

\begin{lemma}
We have $H\vert_{I_n}=P^a(I_n)$.
\end{lemma}

\section{The groups $P$ and $P^\Q$}
In this section, we mainly deal with the group $P$ of all the orientation preserving PL
homeomorphisms of $[0,1]$. In the last part, we remark one word for necessary
modifications with the group $P^\Q$.
For $P$, put $X=(0,1)$ as in Section 4. Then condition
A is trivially fulfilled. Let $H$ be an arbitrary classful subgroup of
$P$. First we shall establish conditions B.

\begin{lemma}\label{l31}
 The group $H$ acts transitively on $(0,1)$.
\end{lemma}

\sc Proof. \rm
There is an element $h_0\in H$ such that $h_0'(0)=2$ and that $h(x)>x$
for any $x\in(0,1)$. Assume $h_0(x)=2x$ on the interval $[0,2^{-n}]$ for
some $n>0$. The interval $[2^{-n-1},2^{-n}]$ is a fundamental domain of the 
action of the group $\langle h_0\rangle$. Thus it suffices to show that
for any 
$\xi\in(2^{-1},1]$, there is an element of $H$ which maps $2^{-n}$ to
$\xi 2^{-n}$. Choose an element $h_1\in H$ such that $h_1'(0)=\xi$.
Assume $h_1$ is linear on an interval $[0,2^{-m}]$ for some $m>n$. Then
$h_1(2^{-m})=\xi2^{-m}$, and hence
$h_0^{m-n}h_1h_0^{n-m}(2^{-n})=\xi2^{-n}$, as is required. \qed

\bigskip
In the rest of this section, we shall establish condition C by the
following lemma.

\begin{lemma}\label{l32}
For some closed interval $I_0\subset(0,1)$, we have $H\vert_{I_0}=P(I_0)$.
\end{lemma}

For any closed interval $I\subset[0,1]$, denote by $\phi_I:[0,1]\to I$ the
orientation preserving bijective linear map.
Define 
$$
P_{1,-1}=\{g\in P\mid g(x)>x, \ \forall x\in(0,1), \ g'(0)=2,\ 
g'(1)=2^{-1}\}.$$
Let $I$ (resp.\ $J$) be a fundamental domain of $g\in P_{1,-1}$ contained in an
end linear zone of $g$ (Definition \ref{elz}) at $0$ (resp.\ at
$1$). Thus 
$I=[a,2a]$ for some $a>0$
and $J=[1-2b,1-b]$ for some $b>0$. If there is $N>0$ such that
$g^N(I)=J$,  we say that the pair $I$ and $J$ are {\em monitoring
intervals for $g$}. The map $f=\phi_J^{-1}g^N\phi_I\in P$ is called the
{\em information of $g$ monitored by $I$ and $J$}. We also say that
{\em $I$ and $J$ monitor the information $f$}.

\begin{definition}
For any $g\in P_{1,-1}$, denote by $\II(g)\subset P$ the set of all the
 monitored
informations of $g$.
\end{definition}

\begin{lemma}\label{l33}
 For any $f\in P$, there is $g\in P_{1,-1}$ such that $f\in \II(g)$.
\end{lemma}

\sc Proof. \rm The proof is almost the same as Lemma \ref{l5}. \qed

\begin{lemma} \label{l34}
Given $g\in P_{1,-1}$ and $f\in \II(g)$, the intervals which
monitor the information $f$ can be chosen arbitrarily near $0$ and $1$. 
\end{lemma}

\sc Proof.\ \rm If the intervals $I$ and $J$ monitor the information
 $f$, and
if $n>0$, then clearly the intervals $g^{-n}(I)$ and $g^n(J)$
 monitor the same information $f$. \qed

\begin{lemma}\label{l35}
 If $g\in P_{1,-1}$ and $g_1\in P$, then $\II(g^{g_1})=\II(g)$.
\end{lemma}

\sc Proof. \rm  Let $g$ and $g_1$ be as in the lemma, and let 
$f\in\II(g)$. It suffices to show that $f\in\II(g^{g_1})$.
By Lemma \ref{l34}, one can choose the monitoring intervals $I$, $J$ of
$g$ which monitor the information $f$ in the end linear zones of $g_1$.
Then $g_1(I)$ and $g_1(J)$ are monitoring intervals of $g^{g_1}$,
with information $f$ since $\phi_{g_1(I)}=g_1\phi_I$ and
$\phi_{g_1(J)}=g_1\phi_J$. \qed

\begin{corollary}\label{c31}
 For any $f\in P$, there is $h\in H\cap P_{1,-1}$ such that $f\in \II(h)$.
\qed
\end{corollary}

Choose an element $h_0\in H\cap P_{1,-1}$ so that ${\rm id}\in \II(h_0)$,
and let $I_0$ and $J_0$ be monitoring intervals of $h_0$ with information
id. That is, there is $N_0>0$ such that
$h_0^{N_0}(I_0)=J_0$ and $\phi_{J_0}^{-1}h_0^{N_0}\phi_{I_0}={\rm id}$.
Let $\hat f$ be an arbitrary element of $P(I_0)$, and let
$f=\phi_{I_0}^{-1}\hat f\phi_{I_0}\in P$. By Corollary \ref{c31},
there is $h_1\in H\cap P_{1,-1}$ such that $f\in \II(h_1)$. Let
$I_1$, $J_1$ be the corresponding monitoring intervals: we assume
$h_1^{N_1}(I_1)=J_1$ for some $N_1>0$ and 
$\phi_{J_1}^{-1}h_1^{N_1}\phi_{I_1}=f$. 
Put
$$ I_0=[a,2a], \ J_0=[1-2b,1-b], \ I_1=[c,2c]\ \mbox{ and }\ J_1=[1-2d,1-d]$$ for some $a,b,c,d>0$.
Choose an element $h_2\in H$ such that $h_2'(0)=c/a$ and $h_2'(1)=d/b$.
Choose a big $n>0$ so that both intervals $h_0^{-n}(I_0)$ and
$h_0^n(J_0)$ are in the end linear zones of $h_2$. Direct computation
shows that
$h_2(h_0^{-n}(I_0))=h_1^{-n}(I_1)$ and
$h_2(h_0^{n}(J_0))=h_1^{n}(J_1)$.
See the figure.
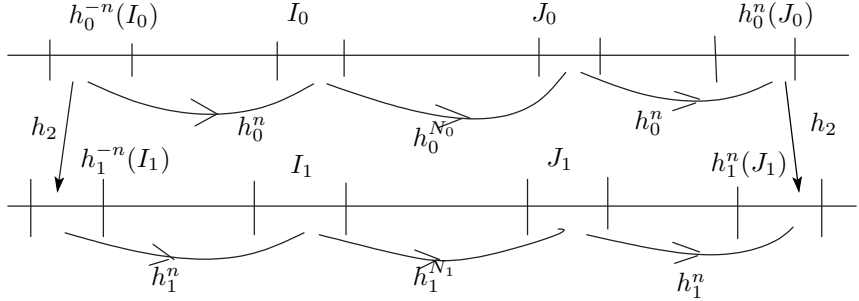
\begin{figure}[h]
{\unitlength 0.1in%
\begin{picture}( 44.4000, 14.3000)(  9.9000,-21.5000)%
%
\special{pn 8}%
\special{pa 990 1010}%
\special{pa 5390 1010}%
\special{fp}%
%
\special{pn 8}%
\special{pa 990 1800}%
\special{pa 5430 1800}%
\special{fp}%
%
\special{pn 8}%
\special{pa 1210 930}%
\special{pa 1210 1140}%
\special{fp}%
%
\special{pn 8}%
\special{pa 1640 950}%
\special{pa 1640 1120}%
\special{fp}%
%
\special{pn 8}%
\special{pa 2400 940}%
\special{pa 2400 1140}%
\special{fp}%
%
\special{pn 8}%
\special{pa 2750 930}%
\special{pa 2750 1120}%
\special{fp}%
%
\special{pn 8}%
\special{pa 3770 920}%
\special{pa 3770 1120}%
\special{fp}%
%
\special{pn 8}%
\special{pa 4090 920}%
\special{pa 4090 1110}%
\special{fp}%
%
\special{pn 8}%
\special{pa 4690 910}%
\special{pa 4700 1150}%
\special{fp}%
%
\special{pn 8}%
\special{pa 5110 920}%
\special{pa 5110 1140}%
\special{fp}%
%
\special{pn 8}%
\special{pa 1110 1660}%
\special{pa 1110 1960}%
\special{fp}%
%
\special{pn 8}%
\special{pa 1490 1660}%
\special{pa 1490 1960}%
\special{fp}%
%
\special{pn 8}%
\special{pa 2280 1670}%
\special{pa 2280 1950}%
\special{fp}%
%
\special{pn 8}%
\special{pa 2760 1680}%
\special{pa 2760 1950}%
\special{fp}%
%
\special{pn 8}%
\special{pa 3710 1660}%
\special{pa 3710 1950}%
\special{fp}%
%
\special{pn 8}%
\special{pa 4130 1650}%
\special{pa 4130 1920}%
\special{fp}%
%
\special{pn 8}%
\special{pa 5250 1660}%
\special{pa 5250 1960}%
\special{fp}%
%
\special{pn 8}%
\special{pa 4810 1700}%
\special{pa 4810 1970}%
\special{fp}%
%
\special{pn 8}%
\special{pa 1410 1150}%
\special{pa 1526 1202}%
\special{pa 1616 1238}%
\special{pa 1646 1249}%
\special{pa 1677 1259}%
\special{pa 1708 1268}%
\special{pa 1740 1277}%
\special{pa 1772 1285}%
\special{pa 1836 1299}%
\special{pa 1869 1304}%
\special{pa 1901 1309}%
\special{pa 1934 1313}%
\special{pa 1967 1316}%
\special{pa 2000 1318}%
\special{pa 2033 1319}%
\special{pa 2066 1319}%
\special{pa 2098 1318}%
\special{pa 2131 1316}%
\special{pa 2163 1313}%
\special{pa 2196 1309}%
\special{pa 2228 1304}%
\special{pa 2259 1298}%
\special{pa 2291 1291}%
\special{pa 2322 1283}%
\special{pa 2382 1263}%
\special{pa 2412 1251}%
\special{pa 2441 1238}%
\special{pa 2499 1210}%
\special{pa 2527 1195}%
\special{pa 2555 1179}%
\special{pa 2583 1164}%
\special{pa 2590 1160}%
\special{fp}%
%
\special{pn 8}%
\special{pa 2660 1160}%
\special{pa 2690 1170}%
\special{pa 2719 1181}%
\special{pa 2839 1221}%
\special{pa 2870 1231}%
\special{pa 2902 1241}%
\special{pa 2933 1251}%
\special{pa 2966 1261}%
\special{pa 2999 1270}%
\special{pa 3032 1280}%
\special{pa 3066 1288}%
\special{pa 3100 1297}%
\special{pa 3134 1305}%
\special{pa 3169 1312}%
\special{pa 3203 1319}%
\special{pa 3238 1325}%
\special{pa 3272 1331}%
\special{pa 3307 1335}%
\special{pa 3341 1339}%
\special{pa 3409 1343}%
\special{pa 3443 1343}%
\special{pa 3476 1342}%
\special{pa 3508 1340}%
\special{pa 3540 1336}%
\special{pa 3572 1331}%
\special{pa 3603 1325}%
\special{pa 3633 1316}%
\special{pa 3662 1306}%
\special{pa 3690 1295}%
\special{pa 3717 1281}%
\special{pa 3744 1266}%
\special{pa 3769 1249}%
\special{pa 3793 1229}%
\special{pa 3816 1208}%
\special{pa 3838 1185}%
\special{pa 3859 1161}%
\special{pa 3901 1111}%
\special{pa 3910 1100}%
\special{fp}%
%
\special{pn 8}%
\special{pa 3990 1100}%
\special{pa 4020 1112}%
\special{pa 4049 1124}%
\special{pa 4079 1135}%
\special{pa 4109 1147}%
\special{pa 4169 1169}%
\special{pa 4229 1189}%
\special{pa 4260 1199}%
\special{pa 4291 1208}%
\special{pa 4322 1216}%
\special{pa 4386 1230}%
\special{pa 4450 1240}%
\special{pa 4483 1244}%
\special{pa 4515 1247}%
\special{pa 4548 1249}%
\special{pa 4580 1250}%
\special{pa 4613 1250}%
\special{pa 4645 1249}%
\special{pa 4678 1246}%
\special{pa 4710 1243}%
\special{pa 4741 1238}%
\special{pa 4773 1231}%
\special{pa 4804 1224}%
\special{pa 4835 1215}%
\special{pa 4865 1205}%
\special{pa 4895 1193}%
\special{pa 4924 1181}%
\special{pa 4953 1167}%
\special{pa 4982 1154}%
\special{pa 4990 1150}%
\special{fp}%
%
\special{pn 8}%
\special{pa 1290 1940}%
\special{pa 1321 1949}%
\special{pa 1351 1958}%
\special{pa 1382 1966}%
\special{pa 1412 1975}%
\special{pa 1443 1983}%
\special{pa 1474 1992}%
\special{pa 1536 2008}%
\special{pa 1567 2015}%
\special{pa 1599 2023}%
\special{pa 1630 2030}%
\special{pa 1662 2037}%
\special{pa 1758 2055}%
\special{pa 1791 2060}%
\special{pa 1823 2064}%
\special{pa 1856 2068}%
\special{pa 1888 2072}%
\special{pa 1921 2074}%
\special{pa 1953 2076}%
\special{pa 1986 2078}%
\special{pa 2050 2078}%
\special{pa 2083 2076}%
\special{pa 2115 2074}%
\special{pa 2146 2071}%
\special{pa 2178 2067}%
\special{pa 2210 2062}%
\special{pa 2272 2048}%
\special{pa 2332 2030}%
\special{pa 2362 2019}%
\special{pa 2422 1995}%
\special{pa 2451 1982}%
\special{pa 2480 1968}%
\special{pa 2509 1955}%
\special{pa 2538 1941}%
\special{pa 2540 1940}%
\special{fp}%
%
\special{pn 8}%
\special{pa 2620 1950}%
\special{pa 2653 1959}%
\special{pa 2687 1967}%
\special{pa 2720 1975}%
\special{pa 2753 1984}%
\special{pa 2785 1992}%
\special{pa 2818 2000}%
\special{pa 2850 2008}%
\special{pa 2912 2024}%
\special{pa 2943 2031}%
\special{pa 2972 2038}%
\special{pa 3002 2045}%
\special{pa 3030 2051}%
\special{pa 3084 2063}%
\special{pa 3110 2068}%
\special{pa 3135 2073}%
\special{pa 3185 2081}%
\special{pa 3210 2083}%
\special{pa 3236 2085}%
\special{pa 3263 2086}%
\special{pa 3291 2086}%
\special{pa 3321 2084}%
\special{pa 3352 2082}%
\special{pa 3385 2078}%
\special{pa 3421 2072}%
\special{pa 3459 2066}%
\special{pa 3500 2057}%
\special{pa 3544 2047}%
\special{pa 3591 2035}%
\special{pa 3642 2021}%
\special{pa 3695 2006}%
\special{pa 3748 1990}%
\special{pa 3797 1974}%
\special{pa 3839 1959}%
\special{pa 3873 1946}%
\special{pa 3895 1934}%
\special{pa 3902 1926}%
\special{pa 3892 1921}%
\special{pa 3870 1920}%
\special{fp}%
%
\special{pn 8}%
\special{pa 4020 1950}%
\special{pa 4144 1978}%
\special{pa 4176 1985}%
\special{pa 4207 1991}%
\special{pa 4238 1998}%
\special{pa 4270 2004}%
\special{pa 4301 2011}%
\special{pa 4365 2023}%
\special{pa 4396 2028}%
\special{pa 4492 2043}%
\special{pa 4525 2047}%
\special{pa 4557 2050}%
\special{pa 4621 2054}%
\special{pa 4653 2055}%
\special{pa 4685 2055}%
\special{pa 4717 2054}%
\special{pa 4749 2052}%
\special{pa 4813 2044}%
\special{pa 4845 2038}%
\special{pa 4877 2031}%
\special{pa 4908 2022}%
\special{pa 4939 2012}%
\special{pa 4969 2000}%
\special{pa 4999 1986}%
\special{pa 5027 1971}%
\special{pa 5053 1952}%
\special{pa 5078 1932}%
\special{pa 5100 1910}%
\special{fp}%
%
\special{pn 8}%
\special{pa 1950 1240}%
\special{pa 2090 1320}%
\special{fp}%
%
\special{pn 8}%
\special{pa 3230 1270}%
\special{pa 3410 1340}%
\special{fp}%
%
\special{pn 8}%
\special{pa 3240 1390}%
\special{pa 3370 1340}%
\special{fp}%
%
\special{pn 8}%
\special{pa 4480 1200}%
\special{pa 4610 1240}%
\special{fp}%
%
\special{pn 8}%
\special{pa 4470 1300}%
\special{pa 4590 1260}%
\special{fp}%
%
\special{pn 8}%
\special{pa 1750 2000}%
\special{pa 1860 2080}%
\special{fp}%
%
\special{pn 8}%
\special{pa 1730 2110}%
\special{pa 1830 2080}%
\special{fp}%
%
\special{pn 8}%
\special{pa 3120 2030}%
\special{pa 3250 2080}%
\special{fp}%
%
\special{pn 8}%
\special{pa 3140 2130}%
\special{pa 3240 2080}%
\special{fp}%
%
\special{pn 8}%
\special{pa 4470 1990}%
\special{pa 4610 2040}%
\special{fp}%
%
\special{pn 8}%
\special{pa 4480 2100}%
\special{pa 4600 2060}%
\special{fp}%
%
\special{pn 8}%
\special{pa 1970 1380}%
\special{pa 2080 1320}%
\special{fp}%
%
\special{pn 8}%
\special{pa 1330 1150}%
\special{pa 1250 1720}%
\special{fp}%
\special{sh 1}%
\special{pa 1250 1720}%
\special{pa 1279 1657}%
\special{pa 1257 1667}%
\special{pa 1239 1651}%
\special{pa 1250 1720}%
\special{fp}%
%
\special{pn 8}%
\special{pa 5060 1150}%
\special{pa 5130 1720}%
\special{fp}%
\special{sh 1}%
\special{pa 5130 1720}%
\special{pa 5142 1651}%
\special{pa 5123 1667}%
\special{pa 5102 1656}%
\special{pa 5130 1720}%
\special{fp}%
\put(13.1000,-8.8000){\makebox(0,0)[lb]{$h_0^{-n}(I_0)$}}%
\put(24.5000,-8.5000){\makebox(0,0)[lb]{$I_0$}}%
\put(37.3000,-8.5000){\makebox(0,0)[lb]{$J_0$}}%
\put(13.7000,-16.3000){\makebox(0,0)[lb]{$h_1^{-n}(I_1)$}}%
\put(24.7000,-16.5000){\makebox(0,0)[lb]{$I_1$}}%
\put(38.1000,-16.1000){\makebox(0,0)[lb]{$J_1$}}%
\put(46.7000,-16.5000){\makebox(0,0)[lb]{$h_1^n(J_1)$}}%
\put(21.9000,-14.6000){\makebox(0,0)[lb]{$h_0^n$}}%
\put(31.1000,-15.2000){\makebox(0,0)[lb]{$h_0^{N_0}$}}%
\put(42.7000,-14.3000){\makebox(0,0)[lb]{$h_0^n$}}%
\put(17.4000,-22.5000){\makebox(0,0)[lb]{$h_1^n$}}%
\put(31.1000,-22.5000){\makebox(0,0)[lb]{$h_1^{N_1}$}}%
\put(44.9000,-22.8000){\makebox(0,0)[lb]{$h_1^n$}}%
\put(11.1000,-14.5000){\makebox(0,0)[lb]{$h_2$}}%
\put(51.9000,-14.3000){\makebox(0,0)[lb]{$h_2$}}%
\put(48.1000,-8.7000){\makebox(0,0)[lb]{$h_0^n(J_0)$}}%
\end{picture}}%

\end{figure}
\\ The equality $\phi_{J_0}^{-1}h_0^{N_0}\phi_{I_0}={\rm id}$
implies (cf.\ the proof of Lemma \ref{l34})
$$\phi_{h_0^{n}(J_0)}^{-1}h_0^{N_0+2n}\phi_{h_0^{-n}(I_0)}={\rm id}.$$
 Likewise we have
$$\phi_{h_1^{n}(J_1)}^{-1}h_0^{N_1+2n}\phi_{h_1^{-n}(I_1)}=f.$$
These equalities show
$$\phi_{h_0^{-n}(I_0)}^{-1}h_0^{-N_0-2n}\phi_{h_0^{n}(J_0)}
\phi_{h_1^{n}(J_1)}^{-1}h_0^{N_1+2n}\phi_{h_1^{-n}(I_1)}=f.
$$
Since $h_2$ is linear on the intervals ${h_0^{-n}(I_0)}$
and ${h_0^n(J_0)}$, we have $h_2\phi_{h_0^{-n}(I_0)}=\phi_{h_1^{-n}(I_1)}$ and
$h_2\phi_{h_0^n(J_0)}=\phi_{h_1^n(J_1)}$. Therefore
$$
\phi_{h_0^{-n}(I_0)}^{-1}h_0^{-N_0-2n}h_2^{-1}h_1^{N_1+2n}h_2\phi_{h_0^{-n}(I_0)}=f.
$$
Finally since $h_0^n\phi_{h_0^{-n}(I_0)}=\phi_{I_0}$, we get

$$
\phi_{I_0}^{-1}h_0^{-N_0-n}h_2^{-1}h_1^{N_1+2n}h_2h_0^{-n}\phi_{I_0}=f.$$
 This implies
$$h_0^{-N_0-n}h_2^{-1}h_1^{N_1+2n}h_2h_0^{-n}\vert_{I_0}=\hat f.
$$
Since $\hat f\in P(I_0)$ is arbitrary, and the map on the LHS is from $H\vert_{I_0}$, the
proof of Lemma \ref{l32} is now complete.
 
\bigskip
For the subgroup $P^\Q$ of $P$ consisting of all the elements with
slopes and breaks in $\Q$, we define $X=\Q\cap(0,1)$. The argument for
$P^\Q$ is the same for $P$ under necessary modifications.


\begin{thebibliography}{99}
\bibitem{Belk} J. M. Belk, {\em Thompson's group $F$}, Thesis, Cornell	University.

\bibitem{BF} {{M. Bestvina and K. Fujiwara, {\em Handlebody subgroups 
        in a mapping class group}, preprint, Arxiv:1412.7847v2.}} 



\bibitem{Gelander} T. Gelander, {\em Convergence groups are not
	invariably generated}, Int.\ Math.\ Res.\ Notices {\bf
	19}(2014), 9806-9814.  
\bibitem{GGJ} T. Gelander, G. Golan and K. Juschenko, {\em Invariable
	generation of Thompson groups}, preprint, Arxiv:1611.08264v2. 
\bibitem{G} R. I. Grigorchuk, {\em On Burnside's problem on periodic
	groups,} Funktional Anal.\ i Prilozhen {\bf 14}(1980), 53-54.  
\bibitem{HNN} G. Higman, B. H. Neumann and H. Neumenn, {\em Embedding
	theorems for groups}, J. London Math.\ Soc.\ {\bf 24}(1949),
	247-254.
\bibitem{KLS} W. M. Kantor, A. Lubotzky and A. Shalev, {\em Invariance
	generations on infinite groups}, J. Algebara {\bf 421}(2015), 296-310.
\bibitem{O} A. Yu.\ Ol'shanskii, {\em Groups of bounded period with
	subgroups of prime order}, Algebra i Logika {\bf 21}(1982), 553-618.
\bibitem{Wiegold1} J. Wiegold, {\em Transformation groups with
	fixed-point-free permutations}, Arch. Math. (Basel) {\bf
	27}(1976), 473-475.
\bibitem{Wiegold2} J. Wiegold, {\em Transformation groups with
	fixed-point-free permutations, II}, Arch. Math. (Basel) 
        {\bf 29}(1977), 571-573.
 



\end{thebibliography}
\end{document}